\documentstyle[fleqn,12pt]{article}
\begin{document}\pagenumbering{arabic}
\setcounter{page}{1}\pagestyle{plain}\baselineskip=16pt
\thispagestyle{empty}

\begin{center}
{\Large\bf 
The Hopf algebra structure of GL(1,H$_q)$ and \\
the isomorphism between SP$_q(1)$ and SU$_q(2)$} 
\end{center}

\vspace{1cm}
Salih Celik\footnote{New E-mail address: sacelik@yildiz.edu.tr}\\
Department of mathematics, Mimar Sinan University, 80690 Besiktas, Istanbul, TURKEY. 

\vspace{1cm}
\begin{center}
{\bf Abstract }
\end{center}
{\footnotesize In this Letter, we introduce the Hopf algebra structure of the 
quantum quaternionic group GL(1,H$_q)$ and discuss the isomorphism 
between the quantum symplectic group SP$_q(1)$ and the quantum unitary group 
SU$_q(2)$.}

\vfill\eject
\noindent
{\bf 1. Introduction }

\noindent
During the past few years, quantum groups [1,2] and $q$-deformed enveloping algebras [3] have been intensively studied both by mathematicians and mathematical physics. From a mathematical point of view, these algebraic structures are just special classes of Hopf algebras [4]. 

In modern mathematical reseach, classical Lie groups and Lie algebras are naturally related not only to the problems of algebras but also to that of geometry and then to that of analysis. These materials are quite basic and play quite important roles in many respects. However, the situation is completely different when we turn our attention to the quantum groups. Quantum groups are not groups in the ordinary sense, they may be thought of as matrix groups in which the matrix elements are themselves noncommutative, obeying sets of bilinear product relations [2]. 

Quaternionic quantum algebra and its coalgebra structure is important in physics. The relevance of quaternionic vector spaces was recently demonstrated in quantum mechanics by Horwitz {\it et al}. [5] and Adler [6]. One of the specific properties of quaternionic mechanics is the dependence of its physical content on the definition of the tensor product [5,7]. 

Quaternionic Lie groups can be defined in terms of matrices 
with quaternionic elements. To make contact with the usual 
formulation of Lie groups in terms matrices with complex matrix 
elements it is useful to represent each quaternion in terms of 
2x2 complex matrices. Therefore, the quantum deformation of the 
quaternion algebra $H$ may be introduce using the idea of the 
quantum matrix theory [1,2].

The notion of quantum quaternions GL$(1,H_q)$ of one parameter 
$q$ was introduced independently by Marchiafava and Rembielinski [8] 
and the author [9]. In this Letter, we obtain the Hopf algebra structure of 
the quantum quaternionic group GL$(1,H_q)$ and show that 
SP$_q(1) \cong $SU$_q(2)$. In Section 2, we introduce the one parameter 
quantum deformation of the quaternionic group GL(1,H) along the 
lines of work of [8] and [9]. In Section 3, we show that there is an 
isomorphism between SP$_q(1)$ and SU$_q(2)$ just as the classical 
case. The Hopf algebra structure of GL$(1,H_q)$ is introduced in 
Section 4.  

\noindent
{\bf 2. Review of GL(1,H$_q)$}

\noindent
We shall define the $q$-deformation of the quaternionic group GL(1,H), 
denoted by GL(1,H$_q$), as an algebra ${\cal A}_q$ equipped with a 
$\star$-operation, where ${\cal A}_q$ has the following properties: 

(1) The unital associative algebra ${\cal A}_q$ generated by $a_0$, $a_1$, 
$a_2$, $a_3$ and $q-$commutation relations [8,9]
$$ x_{\pm} a_k = q^{\pm 1}a_k x_{\pm}, \qquad k=2,3, \qquad 
   a_2a_3 = a_3a_2 \eqno(1) $$
$$ a_0a_1 = a_1a_0 -{{\bf i}\over 2}\lambda_- (a_{2}^{2}+a_{3}^{2}),$$
where ${\bf i} = - 1$, $q$ is a nonzero real number, $\lambda_ = q - q^{-1}$ and $x_{\pm} = a_0 \pm {\bf i} a_1$. 

(2) The operation of conjugation in ${\cal A}_q$ is a map 
$\star :{\cal A}_q \longrightarrow {\cal A}_q$ and acts on the generators 
$a_k$ $(k = 0, 1, 2, ,3)$ of ${\cal A}_q$ as follows [8,9] 
$$ a_{k}^{\star} = a_k, \qquad k=0,1, \qquad 
   a_{2}^{\star} = {1\over 2}(\lambda_{+} a_2 - {\bf i}\lambda_{-} a_3), 
\eqno(2) $$
$$a_{3}^{\star} = {{\bf i}\over 2}(\lambda_{-} a_2 -{\bf i} \lambda_{+} a_3),$$
where $\lambda_{+} = q + q^{-1}$. It is easy to show that the map $\star$ is 
involutive, that is, it satisfies 
$$ (a_{k}^{\star})^{\star} = a_k, \qquad k = 0,1,2,3. \eqno(3) $$
The map $\star$ will be called the anti-involution on ${\cal A}_q$. 

(3) Assume that any quaternion $h$ has a representation 
$$h = a_0 e_0 + a_1 e_1 + a_2 e_2 + a_3 e_3 \eqno(4) $$
with the generators of ${\cal A}_q$ and it will be called the $q$-quaternion, 
and in this case we shall say that $q$-quaternion $h$ belongs to GL(1,H$_q$). 

Here the quaternion units $e_k$ have multiplicative properties defined by 
$$ e_i e_j = - \delta_{ij} e_0 + \epsilon_{ijk} e_k \eqno(5) $$
and 
$$e_0^2 = 1, \qquad e_0 e_k = e_k e_0, \quad k = 1,2,3 \eqno(6)$$
where $\delta_{ij}$ denotes the Kronecker delta and 
$$ \epsilon_{ijk} = {1\over 2}(i - j)(j - k)(k - i) $$
and also 
$$a_k e_l = e_l a_k, \quad k, l = 0,1,2,3. \eqno(7)$$

The quaternionic conjugation is defined by 
$$e_k^\star = 2 \delta_{0,k} e_0 - e_k, \quad k=0,1,2,3. \eqno(8)$$
It is easy to verify that the relations (1) do obey the quantum Yang-Baxter 
equations [1]. 

Relations between starred and unstarred generators are determined by 
$$ x_{\pm} a_{k}^{\star} = q^{\mp 1} a_{k}^{\star} x_{\pm}, \qquad 
   a_{k}^{\star}a_l = a_la_{k}^{\star}, \quad k=2,3, \eqno(9) $$
for $l=2,3$. Note that $(x_\pm)^\star = x_\mp$.

The conjugation of a $q$-quaternion $h \in $GL(1,H$_q$) given by (4) is 
introduced as 
$$\bar{h} = e_0 a_{0}^{\star} - e_1 a_{1}^{\star} - e_2 a_{2}^{\star} - 
  e_3 a_3^\star. \eqno(10)$$
Hence, we can introduce the $q$-{\it norm} of the $q$-quaternion $h$ [8,9] 
$$ {\cal N}_q(h) = h \bar{h} = \bar{h} h = 
  a_0^2 + a_1^2 + {1\over 2} \lambda_+ (a_2^2 + a_3^2) \eqno(11) $$
using the relations (1) and (2) with (5)-(7). It is easy to verify that 
${\cal N}_q(h)$ commutes with all the generators $a_k$ of ${\cal A}_q$, that 
is ${cal N}_q(h)$ belongs to the center of GL(1,H$_q)$.  
Note that if for $h_1$ and $h_2$ are two $q$-quaternions with the generators 
$a_k$ and $b_k$ and the generators $a_k$ commute with $b_k$, then the 
(quaternionic) product $h_1 h_2$ is a $q$-quaternion, i.e., 
$h_1 h_2 \in$ GL(1,H$_q)$. It is also verified that   
$$ {\cal N}_q(h_1 h_2) = {\cal N}_q(h_1) {\cal N}_q(h_2). \eqno(12) $$
Since the $q$-norm, ${\cal N}_q(h)$, is central, the components (the 
generators $a_k$ of ${\cal A}$) of $q-$quaternion $h$ may be normalized so that
${\cal N}_q(h) = 1$, that is, 
$$ a_0 a_1 - a_1 a_0 = {\bf i}{{1-q^2}\over {1+q^2}} (1 - a_0^2 - a_1^2). $$
Thus, we obtain the quantum subgroup SL(1,H$_q)$ = SP$_q(1)$ of quantum 
quaternionic group GL(1,H$_q)$ with ${\cal N}_q(h) = 1$. In the limit 
$q \rightarrow 1$ of the deformation parameter, we obtain a classical 
quaternionic group. 

The addition of two $q$-quaternions must satisfy the following properties in 
order to be a $q$-quaternion. Let $h$ and $h'$ be two $q$-quaternion with the 
generators $a_k$ and $b_k$ for $k=0,1,2,3$ as follows: 
$$x_\pm^1 x_\pm^2 = q^{-2} x_\pm^2 x_\pm^1, \qquad 
  x_+^1 x_-^2  = x_-^2 x_+^1,$$
$$x_+^1 b_k = q^{-1} b_k x_+^1, \qquad 
  x_-^1 b_k = q^{-1} b_k x_-^1 - \lambda_- a_k x_-^2, \quad k=2,3,$$
$$x_-^2 a_k = q a_k x_-^2, \qquad 
  x_+^2 a_k = q a_k x_+^2 + \lambda_- b_k x_+^1, \quad k=2,3,$$
$$x_-^1 x_+^2 = x_+^2 x_-^1 + \lambda_- (y_+^1 y_-^2 + y_+^2 y_-^1),\eqno(13)$$
$$y_+^1 y_+^2 = q^{-2} y_+^2 y_+^1, \qquad 
  y_+^1 y_-^2 = y_-^2 y_+^1 + \lambda_- x_+^1 x_-^2, $$
$$y_-^1 y_-^2 = q^{-2} y_-^2 y_-^1, \qquad 
  y_-^1 y_+^2 = y_+^2 y_-^1  + \lambda_- x_+^1 x_-^2.$$
(It is easy to see that these relations are consistent with the braid statics 
in [10]). Then the generators $c_k = a_k + b_k$ $(k=0,1,2,3)$ for which 
\begin{eqnarray*}
h'' 
& = & h + h' = (x_+^1 + x_+^2) + (y_+^1 + y_+^2) {\bf j} \\
& = & (a_0 + b_0) e_0 + (a_1 + b_1) e_1 + (a_2 + b_2) e_2 + (a_3 + b_3) e_3 
\end{eqnarray*}
also obey the relations (1) so that $h''$ is in GL(1,H$_q$). Here 
${\bf j} = e_2 ~({\bf i j} = e_3)$. 

\noindent
{\bf 3. The Isomorphism Between SP$_q(1)$ and SU$_q(2)$ }

\noindent
We begin with some information. It is well known that if 
$\varphi : G_1 \longrightarrow G_2$ is a one-to-one homomorphism of groups, 
then the map $\varphi$ is an isomorphism of $G_1$ onto the subgroup 
$\varphi(G_2)$ of $G_2$. So we can consider $G_1$ as a subgroup of $G_2$. 
First we will consruct a one-to-one homomorphism 
$$\varphi : \mbox{GL(1,H}_q) \longrightarrow \mbox{GL}_q(2) \eqno(14)$$
and then for $h \in$ GL(1,H$_q$) we are going to use as the $q$-norm of $h$ 
the $q$-determinant of $\varphi(h)$. Here the group GL$_q(2)$ is the quantum 
group of 2x2 nonsingular matrices whose matrix elements obey certain 
$q$-dependent commutation relations (similar to (1)) [2]. 

For each $q-$quaternion $h$ given by (4), define $\varphi(h)$ as a 2x2 matrix 
by 
$$ \varphi(h) = \left ( \matrix{ a_0 + {\bf i}a_1 & a_2 + {\bf i}a_3 \cr \cr 
                             -a_2 + {\bf i}a_3 & a_0 - {\bf i}a_1 \cr}
\right), \eqno(15)$$
where 
$$\varphi(e_1) = {\bf i} \sigma_3, \qquad 
  \varphi(e_2) = {\bf i} \sigma_2, \qquad 
  \varphi(e_3) = {\bf i} \sigma_1$$
with the Pauli matrices satisfying 
$$ \sigma_i \sigma_j = \delta_{ij} + {\bf i} \epsilon_{ijk} \sigma_k, \qquad 
   \sigma_i \sigma_j + \sigma_j \sigma_i = 0, \quad i \neq j. $$
It is easy to verify that the map $h \longmapsto \varphi(h)$ is 
one-to-one and
$$ \varphi(h_1 h_2) = \varphi(h_1) \varphi(h_2) \eqno(16) $$
for $h_1, h_2 \in$ GL(1,H$_q)$ provided that the components (the generators of 
${\cal A}_q$) of $h_1$ commute with those of $h_2$ (recall that, in this case 
$h_1 h_2 \in$ GL(1,H$_q)$). So this map is a homomorphism. 

Next, if $h \in$ GL(1,H$_q$) we have 
$$Det_q(\varphi(h)) = (a_0 + {\bf i} a_1) (a_0 - {\bf i} a_1) + 
    q (a_2 + {\bf i} a_3) (a_2 - {\bf i} a_3) = {\cal N}_q(h). \eqno(17)$$

Now we discuss the isomorphism between SP$_q(1)$ and SU$_q(2)$. SP$_q(1)$ is 
the set of all $q$-quaternions of unit length ($q$-norm) and SU$_q(2)$ is the 
set of all complex 2x2 matrices $A$ such that (see, for example, [11]) 
$$ (A^{\star})^T A = I \qquad \mbox{and} \qquad Det_q(A) = 1, \eqno(18) $$
where T denotes the matrix transposition and $A^{\star}$ stands for the 
Hermitean conjugate of $A$. That is, $A \in$ SU$_q(2)$ if and only if 
$A \in$ GL$_q(2)$ with (18). The operation in SP$_q(1)$ is multiplication of 
$q$-quaternions, in SU$_q(2)$ it is matrix multiplication. Consider the map
$$ \varphi : \mbox{SP}_q(1) \longrightarrow \mbox{SU}_q(2). \eqno(19) $$

We have seen that the map $\varphi$ induces a one-to-one homomorphism of 
GL(1,H$_q$) into GL$_q(2)$, thus restriction of $\varphi$ to SP$_q(1)$ is 
still a ono-to-one homomorphism. Therefore, we just need to show that 

{\bf (1)}
$h \in \mbox{SP}_q(1) ~~\Longrightarrow ~~ \varphi(h) \in \mbox{SU}_q(2)$ and

{\bf (2)}
for every $A \in \mbox{SU}_q(2)$ there exists some $\varphi(h)$ with 
$h \in \mbox{SP}_q(1)$.

Note that the matrix $\varphi(h)$ in (15) is a GL$_q(2)$ matrix. This is 
easy to show that using the relations (1). 

Let 
$$ A = \left(\matrix{ a & b \cr c & d \cr} \right) $$
be in SU$_q(2)$. Then we find that 
$$d = a^\star \qquad \mbox{and} \qquad c = - q^{-1} b^\star \eqno(20)$$
with Det$_q(A) = 1$. 

On the other hand, using the relations (2), we can write 
$([\varphi(h)]^{\star})^T \varphi(h) = I$ since ${\cal N}_q(h) = 1$. Here $I$ 
denotes the 2x2 unit matrix. Also Det$_q(\varphi(h)) = 1$. Now 
$$ a_2 - {\bf i} a_3 = q^{-1} (a_2 + {\bf i} a_3)^{\star} \eqno(21) $$
with (2), so that
$$ \varphi(h) = \left ( \matrix{ a_0 + {\bf i}a_1 & a_2 + {\bf i}a_3 \cr \cr 
             -q^{-1} (a_2 + {\bf i}a_3)^{\star} & (a_0 + {\bf i}a_1)^{\star} 
\cr} \right). \eqno(22) $$
It is easy to verify that the matrix elements of $\varphi(h)$ satisfy the 
relations of SU$_q(2)$ provided $(ab)^{\star} = b^{\star}a^{\star},$ that is, 
$\varphi(h) \in \mbox{SU}_q(2)$. So, if $a = a_0 + {\bf i} a_1$ and 
$b = a_2 + {\bf i} a_3$ we may take 
$h = a_0 + {\bf i} a_1 + (a_2 + {\bf i} a_3) {\bf i}$ and have 
$\varphi(h) = A$ (and ${\cal N}_q(h) = 1$). This proves that the map $\varphi$ 
in (19) is an isomorphism. 

\noindent
{\bf 4. The Hopf Algebra Structure of $GL(1,H_q)$ }

\noindent
Now we introduce three operators $\Delta$, $\epsilon$ and $\cal S$ on 
GL(1,H$_q)$, which are called the comultiplication, the counit, and the 
antipode (co-inverse), respectively. 

{\bf (1)}
The comultiplication $\Delta$ is defined by
$$ \Delta(h) = h \otimes h. \eqno(23) $$
The action of comultiplication $\Delta$ on the generators $a_k$ of $\cal A_4$ 
can be introduced as follows:
$$ \Delta(a_0) = a_0 \otimes a_0 - (a_1 \otimes a_1 + a_2 \otimes a_2 + 
a_3 \otimes a_3), $$
$$ \Delta(a_1) = a_0 \otimes a_1 + a_1 \otimes a_0 + a_2 \otimes a_3 - 
a_3 \otimes a_2 , \eqno(4.2) $$
$$ \Delta(a_2) = a_0 \otimes a_2 + a_2 \otimes a_0 + a_3 \otimes a_1 - 
a_1 \otimes a_3 , $$
$$ \Delta(a_3) = a_0 \otimes a_3 + a_3 \otimes a_0 + a_1 \otimes a_2 - 
a_2 \otimes a_1, $$
and 
$$ \Delta(e_0) = e_0 \otimes e_0, \eqno(4.3) $$
where $\otimes$ denotes the tensor product. Note that the relations (24) are 
invariant under (1). For simplicity, we prove the invariance of one of the 
relations in (1) here. Proofs of the remaining formulas are similar. A 
direct calculation shows that 
\begin{eqnarray*}
\Delta(x_+)\Delta(a_2) 
& = & q \Delta(a_2) (x_+ \otimes x_+) 
    - {\lambda_- \over 2} [\Delta(a_2) + {\bf i} \Delta_2] (y_+ \otimes y_-) \\
&  & - {1\over 2} [\lambda_+ \Delta(a_2) - {\bf i} \lambda_- \Delta_2] 
     (y_+ \otimes y_-) \\
& = & q \Delta(a_2) \Delta(x_+) 
\end{eqnarray*} 
where $y_\pm = a_2 \pm {\bf i} a_3$ and 
$$\Delta_2 = a_0 \otimes a_3 - a_3 \otimes a_0 + a_1 \otimes a_2 + 
  a_2 \otimes a_1. $$

Using (23), it is also easy to show that 
$$ \Delta({\cal N}_q(h)) = {\cal N}_q(h) \otimes {\cal N}_q(h). \eqno(26) $$
The comultiplication $\Delta$ is an algebra homomorphism which is 
co-associative, that is, 
$$ ({\cal I} \otimes \Delta) \circ \Delta = 
  (\Delta \otimes {\cal I}) \circ \Delta, 
\eqno(4.5) $$
where $\circ$ stands for the composition of maps and  
${\cal I}$ is the identity map. 

{\bf (2)} 
The counit $\varepsilon$ is introduced by
\footnote{The relations (23) and (28) also appear in the paper of 
Marchiafava and Rembielinski [8]}
$$ \varepsilon(h) = e_0 \eqno(28) $$
whose action on the generators $a_k$ of $\cal A_q$ can be defined by
$$ \varepsilon(a_k) = \delta_{0,k} e_0, \quad k=0,1,2,3 \eqno(29) $$
and also 
$$ \varepsilon(e_0) = e_0. \eqno(30)$$
The counit $\varepsilon$ is an algebra homomorphism such that
$$ (\varepsilon \otimes {\cal I}) \circ \Delta = {\cal I} = ({\cal I} \otimes 
  \varepsilon) \circ \Delta. $$
Thus we have verified that GL(1,H$_q)$ is a {\it bialgebra} with 
multiplication $m$ satisfying the associativity axiom:
$$ m \circ (m \otimes {\cal I}) = m \circ ({\cal I} \otimes m) $$
where $m(h_1 \otimes h_2) = h_1 h_2$.

A bialgebra with the extra structure of the antipode is called a Hopf 
algebra [4].

{\bf (3)}
We can introduce the antipode ${\cal S}$ as follows: 
$$ {\cal S}(h) = h^{-1}. \eqno(31) $$
The action of antipode ${\cal S}$ on the generators $a_k$ of $\cal A_q$ can be 
defined by
$$ {\cal S}(a_k) = {\cal N}_{q}^{-1}(h) (2 \delta_{0,k}a_0 - a_{k}^{\star}) 
\eqno(32) $$
for $k = 0,1,2,3$. The antipode ${\cal S}$ is an algebra anti-homomorphism 
which satisfies
$$ m \circ [({\cal S} \otimes {\cal I}) \circ \Delta] = \varepsilon = 
  m \circ [({\cal I} \otimes {\cal S}) \circ \Delta]. \eqno(33) $$

The comultiplication, counit and antipode which are specified above supply 
GL(1,H$_q)$ with a Hopf algebra structure. 

\noindent
{\bf 5. Discussion }

\noindent
We have introduced the quantum deformation of one parameter of 1x1 
quaternionic group, GL(1,H$_q)$ along the lines of the work of [9]. 
We have constructed the Hopf algebra structure of GL(1,H$_q)$ and have 
discussed an isomorphism between SP$_q(1)$ and SU$_q(2)$. We hope that the 
methods used in this paper will be helpful for an explicit construction of 
GL(n,H$_q)$ and its Hopf algebra structure. 

\noindent
{\bf Acknowledgement}

\noindent
This work was supported in part by TBTAK the Turkish Scientific and Technical 
Research Council. 

\noindent
{\bf References }

\begin{description}
\item [1.] Reshetikhin, N. Y., Takhtadzhyan, L. and Faddeev, L.: 
  {\it Leningrad Math. J} {\bf 1} (1990), 193-225. 
\item[2.] Manin, Yu I.: {\it Quantum groups and noncommutative geometry}, 
  Montreal Univ. Preprint, 1988.
\item[3.] Drinfeld, V.: in {\it Proc. Int. Cong. Math. Berkeley}, 1986. 
\item[4.] Sweetler, M. E.: {\it Hopf Algebras}, Benjamin, New York, 1969. 
\item[5.] Horwitz, L. P. and Biedenharn, L. C.:{\it Ann. Phys.} {\bf 157} 
  (1984), 432-488.
\item[6.] Adler, S. L.: {\it Comm. Math. Phys.} {\bf 104} (1986), 611-656.
\item[7.] Nash, C. G. and Joshi, G. C.: {\it J. Math. Phys.} {\bf 28} (1987), 
  2883-2885.
\item[8.] Marchiafava, S. and Rembielinski, J.: {\it J. Math. Phys.} 
  {\bf 33} (1992), 171-173. 
\item[9.] Celik, S.: {\it The quantum matrix groups and $q$-oscillators} 
  (in Turkish), PhD Thesis, Istanbul Technical University, 1992. 
\item[10.] Majid, S.: {\it J. Math. Phys.} {\bf 35} (1994), 2617-2632.
\item[11.] Woronowicz, S. L.: {\it Publ. RIMS Kyoto Univ.} 23 (1987), 117-181.

\end{description}

\end{document}